\documentclass{birkjour}
 \newtheorem{thm}{Theorem}[section]

 \theoremstyle{definition}
 
 \theoremstyle{remark}
 \newtheorem{rem}[thm]{Remark}
 
 \numberwithin{equation}{section}

\usepackage{amsmath, amsthm, amssymb}
\usepackage{url} 
\usepackage{graphicx}
\usepackage{xcolor}

\def\r{\mathbb R}
 
\def\s{\mathbb S}
\def\h{\mathbb H}
\def\m{\mathbb M}
\def\t{\mathbf t}

\begin{document}

\title[Curves  on
a Totally Umbilical Surface of a Space Form]
 {A Characterization of Curves that Lie on
a Totally Umbilical Surface of a Space Form}

\author[Rafael L\'opez]{Rafael L\'opez}
\address{Departamento de Geometr\'{\i}a y Topolog\'{\i}a Universidad de Granada 18071 Granada, Spain}
\email{rcamino@ugr.es}


\begin{abstract}  We give necessary and sufficient conditions on the curvature and the torsion of a regular curve of the space forms $\h^3$ and $\s^3$ to be contained in a   totally umbilical surface.  In case that the curve has  constant torsion, we obtain the value of  the curvature of the curve. Also numerical pictures of these curves are shown.
 
\end{abstract}

\maketitle

\section{Introduction and statement of the results} \label{intro}

In this paper, we investigate the following  problem: 
\begin{quote} Given a surface $S$ in a $3$-dimensional space $M^3$, how to characterize   a regular curve $\alpha\colon I\subset\r\to M^3$ to be contained in $S$ in terms of the curvature and the torsion of $\alpha$? 
\end{quote}
Although the formulation of this problem is simple, surprisingly  its solution, if any, is  hard and only in a few cases of surfaces $S$  a satisfactory answer   is known.  A regular curve $\alpha\colon I\subset\r\to M^3$, which it is not a geodesic, has defined its curvature $\kappa$ and its torsion $\tau$. If the ambient space $M^3$ has constant curvature, both   quantities $\kappa$ and $\tau$ determine the curve $\alpha$ up to a rigid motion of $M^3$. For general manifolds $M^3$ the uniqueness modulo congruences is valid once is prescribed initial conditions on the position and the Frenet frame \cite{cas0,cas1}. 
 
Let now $S$ be a given surface in $M^3$. If the curve $\alpha$ is contained in $S$, that is, $\alpha(I)\subset S$, there is not an {\it a priori} condition on $\kappa$ and $\tau$ that characterizes the fact that the image of $\alpha$ lies contained in $S$ because the notions of $\kappa$ and $\tau$ are intrinsic of the curve.  

An interesting application of the above problem is in computer vision  \cite{fo,ij,jt,krsw,kp}.   Suppose that we have a spatial curve in Euclidean space of which we have sufficient information data with respect to the position of its points. From these data, we are able to compute the curvature and the torsion of the curve. The problems asks whether such a curve is contained on a given  surface $S$ of Euclidean space.  As one can notice, the problem is difficult because the data corresponds to a one-dimensional object, whereas the object we are looking for is two-dimensional.

 Without loss of generality, we will assume in this paper  that all curves are  parametrized by arc length, being $s$ is the arc length parameter. We also assume sufficiently regularity of the curves, at least $C^4$, although one may demand something weaker that $C^4$ by choosing an appropriate frame different than the Frenet frame. 
 
Following with the Euclidean space $\r^3$, two cases of surfaces $S$ have a known answer to the problem:
\begin{enumerate}
\item   {\it $S$ is a plane}. Then $\alpha$ lies contained in $S$ if and only if $\tau=0$.
\item   {\it $S$ is a sphere}. Suppose that $R$ is the radius of the sphere. If $\alpha$ lies contained in $S$ with $\tau\not=0$, then  
\begin{equation}\label{es}
\frac{1}{\kappa^2}+\frac{\kappa'^2}{\kappa^4\tau^2}=R^2,
\end{equation}
where the prime $(')$ stands for the derivative with respect to  $s$. If $\kappa'\not=0$, then the converse is true. If the radius of the sphere is unknown, then   $\alpha$ lies contained in some sphere of $\r^3$ if and only if 
\begin{equation}\label{es2}
\frac{\tau}{\kappa}-\left(\frac{\kappa'}{\kappa^2\tau}\right)'=0.
\end{equation}
\end{enumerate}
 Notice that Eq. \eqref{es} is a     third order differential equation in the coordinates of the curve whereas \eqref{es2} is of fourth order. Equations \eqref{es} and \eqref{es2} appear in many textbooks of differential geometry: see for example, Ex.  13, Sect.  1.5, in \cite{dc}.  Since there is more than one way to prove   \eqref{es} and \eqref{es2}, we review different approaches given in the literature. For example, in \cite[p. 32]{st}, Struik exploits the idea of contact  between curves and surfaces to obtain \eqref{es} but to prove \eqref{es2} Struik assumes that $\kappa'\not=0$ as in \cite{dc}.   K\"{u}hnel in \cite[Th. 2. 10]{ku} also uses the idea of contact to characterize  with \eqref{es2} spherical curves by using that if a curve $\alpha$ has a contact of order $k$ with a surface $S=\{{\bf x}\in\r^3:G({\bf x})=0\}$ then  all derivatives of $s\mapsto  G(\alpha(s))$ up to order $k$ vanish (see also \cite[Lem. 19.2]{kr}).  Yet another proof using order of contact is given by da Silva \cite[Th. 2.1]{ds20} where the author applies the theory of contact as K\"{u}hnel did but using Bishop frames. Other question related with  \eqref{es} and \eqref{es2}  is that one assumes that the torsion $\tau$ does not have zeros. A characterization  that does not require the torsion to be non-zero can be found in \cite{wo}.  Extensions of \eqref{es} and \eqref{es2} when the ambient space is the  Lorentz-Minkowski space $\mathbb{L}^3$ and the surface $S$ is the pseudo-sphere or the hyperbolic plane  appear in \cite{da,il,pp,pm1,pm2,pm3}.

Other simple surface $S$ where we can pose the problem in Euclidean space is a circular cylinder, $S=\{(x,y,z)\in\r^3:x^2+y^2=R^2\}$, $R>0$. However, the problem is far to have a successful answer. In \cite{krsw}  the authors give   necessary conditions for $\alpha$ to be included in $S$ in terms of a resultant of two polynomials whose coefficients  are expressed by $\kappa$, $\tau$ and their first and second derivatives. More recently,  in \cite{sh} the authors derive some conditions when  $\alpha$ is regarded  as a curve at constant separation from the axis of the cylinder. Then they obtain a differential equation involving $\kappa$, $\tau$ and the angle $\xi$ between the separation vector from the axis and the normal vector  of $\alpha$: see Props. 4 and 5 in \cite{sh}.

A different approach to the initial problem is using different notions of curvatures.  Bishop, in a remarkable paper \cite{bi}, introduced a new frame for a spatial curve $\alpha\colon I\subset\r\to\r^3$. In the new frame $\{T=\alpha',N_1,N_2\}$ the key property is that $T$ and $N_i'$ are parallel, $i=1,2$. In the corresponding Bishop equations, two new curvatures  $\kappa_1$ and $\kappa_2$ are defined.  Comparing with Frenet frames, Bishop frames exist even if the curve has points with $\kappa=0$.   Bishop characterized curves contained in spheres by means of a linear relation of $\kappa_1$ and $\kappa_2$.  Generalizations of this result have been obtained when $S$ is a geodesic sphere  of $\s^3$ and $\h^3$  \cite{dd} or a pseudosphere of   the Lorentz-Minkowski space   $\mathbb{L}^3$ \cite{da}.  

In this paper we consider the initial problem when $M^3$ is the hyperbolic space $\h^3$ or the $3$-sphere $\s^3$ and $S$ is a totally umbilical surface. More precisely, we extend the characterizations \eqref{es} and \eqref{es2}  for all totally umbilical surfaces of $\h^3$ and $\s^3$. To distinguish  with the use of Bishop frames given in \cite{da,dd}, we will follow the same ideas of \cite{dc} using Frenet frames and requiring that $\tau$ has not zeros. Notice that $\tau\not=0$ is an open condition so if $\tau(s_0)\not=0$ at some point $s_0$, then $\tau\not=0$ in a subinterval of $\alpha$ around $s_0$. 

Let  $\m^3(c)$ be the space form of constant sectional curvature $c \in\{-1,0,1\}$.   Let $S\subset\m^3(c)$ be a totally umbilical surface. If $c=0$, then $S$ is  a plane or a round sphere and the problem has been solved as it has been previously explained. If $c=\pm 1$ and when $S$ is a geodesic sphere, equation \eqref{es2} appeared  in \cite[Thm. 4]{dd} (see Rem. \ref{r2} below). In the present paper and for the use of Frenet frames, we use the hyperquadric model of $\h^3$ and $\s^3$ and Eq. \eqref{es} will be extended for all totally umbilical surfaces. The value of $R$ will be substituted in terms of the mean curvature $H$ of $S$.

We consider the hyperquadric  model of $\h^3$ and  $\s^3$.    In   the vector space $\r^4$ with standard coordinates $p=(x_1,x_2,x_3,x_4)$, define the metric 
\begin{equation}\label{me}
\langle,\rangle= dx_1^2+dx_2^2+dx_3^2+c\, dx_4^2.
\end{equation}
Consider the hyperquadric $\mathcal{H}=\{p\in \r^4:\langle p,p\rangle=c\}$. In case that  $c=-1$, we only take the component of $\mathcal{H}$  where the $x_4$-coordinate is positive. Then $\m^3(c)$ is isometric to $\mathcal{H}$ endowed with the induced metric $\langle,\rangle$ from $ \r^4$.  Notice that if $c=1$, then $\mathcal{H}=\s^3$ is the round sphere in Euclidean space $\r^4$ of radius $1$ and centered at the origin.

The totally umbilical surfaces of $\m^3(c)$ are the intersection of affine hyperplanes of $\r^4$ with $\mathcal{H}$, that is, 
$$U_{\vec{a},\sigma}=\{p\in\mathcal{H}:\langle p,\vec{a}\rangle=\sigma\},$$
where $\sigma\in\r$ and $\vec{a}\in \r^4$, $\vec{a}\not=0$, with $\langle\vec{a},\vec{a}\rangle=\epsilon\in\{-1,0,1\}$ \cite{sp}.  If  $c=1$, then necessarily $\epsilon=1$ and $\sigma\in (-1,1)$. We give the  classification of the totally umbilical surfaces. 
For distinguishing if necessary, we add a superscript $()^h$ and $()^s$ to the notation of $U_{\vec{a},\sigma}$ depending if the surface is contained in $\h^3$ or $\s^3$, respectively.

\begin{figure}[hbtp]
\centering
\includegraphics[width=.6\textwidth]{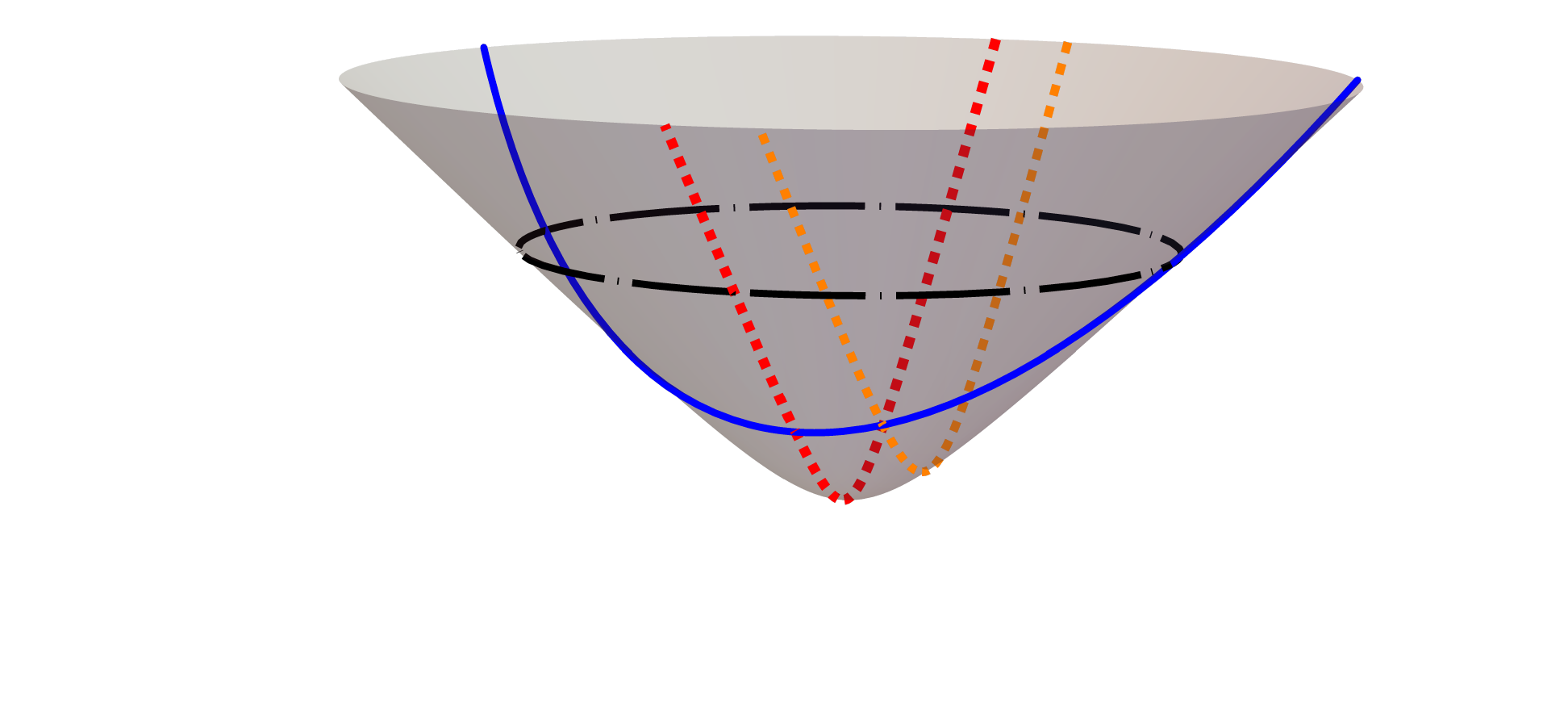}\includegraphics[width=.35\textwidth]{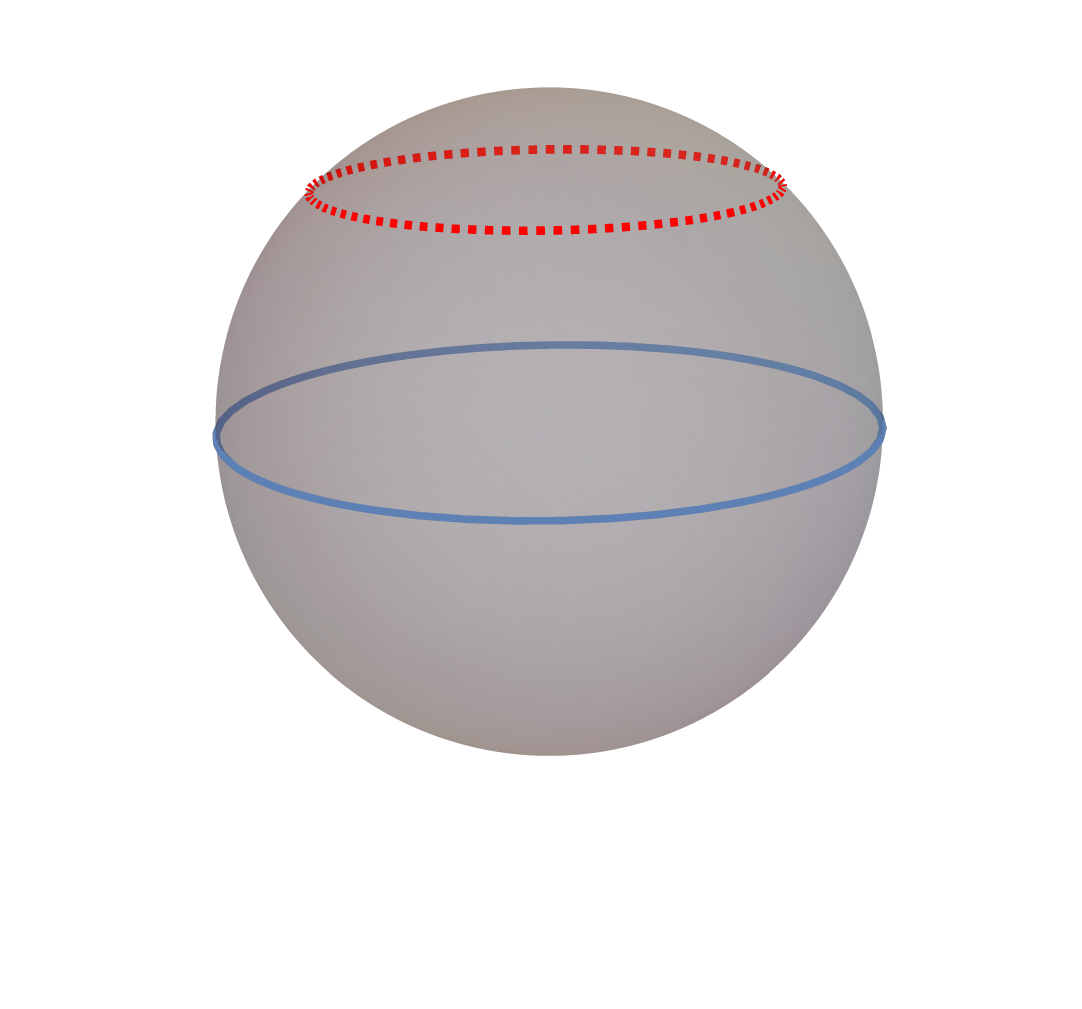}
\caption{Left: Totally umbilical surfaces in the hyperbolic space $\h^3$: totally geodesic planes and equidistant surfaces (red and orange dotted curves), horospheres (blue solid curve) and geodesic spheres (black dashed curve). Right: totally umbilical surfaces in the sphere $\s^3$: totally geodesic spheres (blue solid curve) and geodesic spheres (red dotted curve).   }\label{fig01}
\end{figure}

Consider the hyperbolic space $\h^3$. In this case $c=-1$ and  the metric   \eqref{me} is Lorentzian. The classification of the totally umbilical surfaces depends on the causal character of the affine hyperplane:
 \begin{enumerate}
 \item totally geodesic planes: $\epsilon=1$, $\sigma=0$;
 \item equidistant surfaces: $\epsilon=1$, $\sigma\not=0$;
 \item horospheres: $\epsilon=0$, $\sigma\not=0$;
 \item geodesic spheres: $\epsilon=-1$, $|\sigma|>1$.
 \end{enumerate}
 See Fig. \ref{fig01}, left. A unit normal vector field on $U_{\vec{a},\sigma}^h$ is 
 $$\xi(p)=-\lambda(\vec{a}+\sigma p),\quad \lambda=\frac{1}{\sqrt{\epsilon+\sigma^2}}.$$
 If $A$ is the shape operator on $U_{\vec{a},\sigma}^h$, then $A\xi=\lambda\sigma\cdot{\rm id}$, which proves that $U^h_{\vec{a},\sigma}$ is totally umbilical. Notice that the mean curvature of $U^h_{\vec{a},\sigma}$ is \begin{equation}\label{hh}
 H=\lambda\sigma=\frac{\sigma}{\sqrt{\epsilon+\sigma^2}}
 \end{equation}
 and the extrinsic curvature is $K_{ex}=\sigma^2/(\epsilon+\sigma^2)$. The Gauss equation implies that the intrinsic curvature $K$ is constant, namely, $K=K_{ext}-1=-\frac{\epsilon}{\epsilon+\sigma^2}$.
 
 \begin{thm}\label{t1} Let $U_{\vec{a},\sigma}^h$ be a totally umbilical surface of $\h^3$. Let  $\alpha\colon I\subset\r\to\h^3$   be a regular curve.   
 \begin{enumerate}
 \item Case $\sigma=0$. If $\alpha$ is contained  in  $U_{\vec{a},0}^h$ then  $\tau=0$.
  \item Case $\sigma\not=0$. If $\alpha$ is contained in   $U_{\vec{a},\sigma}^h$ with $\tau\not=0$, then
 \begin{equation}\label{eq1}
\frac{1}{\kappa^2}+\frac{\kappa'^2}{\kappa^4\tau^2}=\frac{1}{H^2},
 \end{equation}
 \end{enumerate}
 where $H$ is the mean curvature of $U_{\vec{a},\sigma}^h$.  
 
 Reciprocally, if $\alpha$ has constant zero torsion $\tau=0$ then $\alpha$ is contained in a totally geodesic plane. In case that $\tau\not=0$ and $\kappa'\not=0$, if  
 \begin{equation}\label{eq13}
\frac{1}{\kappa^2}+\frac{\kappa'^2}{\kappa^4\tau^2}=C>0 
 \end{equation}
 for some constant $C$, then  $\alpha$ is contained in a totally umbilical surface of mean curvature $H$, where $C=1/H^2$:  an equidistant surface if $C>1$, a horosphere if $C=1$ or a geodesic sphere if $C<1$.  \end{thm}

If it is unknown the type of totally umbilical surface, the  analogous result to \eqref{es2} is the following.

 \begin{thm}\label{t12}  Let   $\alpha\colon I\subset\r\to\h^3$   be a regular curve. Then $\alpha$ is contained in a totally geodesic plane if and only if $\tau=0$. If $\tau\not=0$ and $\kappa'\not=0$, then $\alpha$ is contained in a totally non-geodesic umbilical surface  if and only if  
 \begin{equation}\label{eq12}
\frac{\tau}{\kappa}-\left(\frac{\kappa'}{\kappa^2\tau}\right)'=0.
 \end{equation}
 \end{thm}
 
 Consider now the unit sphere $\s^3$. The totally umbilical surfaces $U^s_{\vec{a},\sigma}$ are totally geodesic spheres  (or great spheres) if $\sigma=0$ or geodesic spheres if $\sigma\not=0$. In the latter case, the radius of the sphere is   $R=\sqrt{1-\sigma^2}$. See Fig. \ref{fig01}, right.  Now the unit normal vector of $U^s_{\vec{a},\sigma}$ is $\xi(p)=(a-\sigma p)/\sqrt{1-\sigma^2}$ and $A\xi=\frac{\sigma}{\sqrt{1-\sigma^2}}\cdot {\rm id}$. The mean curvature is $H=\frac{\sigma}{\sqrt{1-\sigma^2}}$ and the intrinsic curvature   is $K=\frac{1}{1-\sigma^2}$.
 
  \begin{thm}\label{t2} Let $U_{\vec{a},\sigma}^s$ be a totally umbilical surface of $\s^3$. Let $\alpha\colon I\subset\r\to\s^3$   be a regular curve. 
   \begin{enumerate}
 \item Case $\sigma=0$. If $\alpha$ is contained  in  $U_{\vec{a},0}^s$ then $\tau=0$.
  \item Case $\sigma\not=0$. If $\alpha$ is contained  in  $U_{\vec{a},\sigma}^s$ with $\tau\not=0$, then 
   \begin{equation}\label{eq2}
\frac{1}{\kappa^2}+\frac{\kappa'^2}{\kappa^4\tau^2}=\frac{1}{H^2},
 \end{equation}
 \end{enumerate}
 where $H$ is the mean curvature of $U_{\vec{a},\sigma}^s$. 
 
 Reciprocally, if $\alpha$ has constant zero torsion $\tau=0$ then $\alpha$ is contained in a totally geodesic sphere. In case that $\tau\not=0$ and $\kappa'\not=0$, if  
  \begin{equation}\label{eq14}
\frac{1}{\kappa^2}+\frac{\kappa'^2}{\kappa^4\tau^2}=C>0 
 \end{equation}
 for some constant $C$, then  $\alpha$ is contained in some geodesic sphere of mean curvature $H$ with $C=1/H^2$.
 \end{thm}

 \begin{thm}\cite{dd} \label{t22}Let $\alpha\colon I\subset\r\to\s^3$   be a regular curve. Then  $\alpha$ is contained in a totally geodesic sphere if and only if $\tau=0$. If $\tau\not=0$ and $\kappa'\not=0$, then $\alpha$ is contained in a geodesic sphere if and only if  
 \begin{equation}\label{eq3}
\frac{\tau}{\kappa}-\left(\frac{\kappa'}{\kappa^2\tau}\right)'=0
 \end{equation}
 \end{thm}
 
 In all these results, the dimension of the ambient space is $3$. For higher dimensions and in arbitrary manifolds $M^{n+1}$,    
   the generalization of the Bishop frame, the so-called rotation minimizing frame (RM frame in short), also is useful in characterizations of curves. A RM frame along a curve $\alpha$ is a frame $\{T,N_1,\ldots,N_n\}$ where $T$ and $\nabla_TN_i$ are parallel.  In \cite[Thm. 3]{dd}, the authors prove that a regular curve $\alpha$ lies on a geodesic hypersphere of radius $R$ of $\h^{n+1}$ (resp. $\s^{n+1}$) if and only if $\sum_{i=1}^na_i\kappa_i+\coth(R)=0$ (resp. $\sum_{i=1}^na_i\kappa_i+\cot(R)=0$). Here $a_i$ are real constants and $\kappa_i$ are the curvatures of the RM frame. In the present paper, and using the same ideas,  we extend this result to curves contained in totally umbilical hypersurfaces of $\m^{n+1}(c)$, $c\not=0$ (for geodesic hyperspheres, it is the same result that \cite{dd}). Notice that now the regularity that we demand is weaker than $C^4$, and $C^2$   is enough. 

\begin{thm}\label{t5}
Let  $\alpha\colon I\subset\r\to\m^{n+1}(c)$ be a $C^2$ curve, $c\in\{-1,1\}$. If   $\alpha$ is contained in a totally umbilical hypersurface $U_{\vec{a},\sigma}$, then   there is a linear relation
\begin{equation}\label{eq-rm}
\sum_{i=1}^na_i\kappa_i+\sigma=0,
\end{equation}
where $a_i\in\r$, $1\leq i\leq n$, are real constants. Reciprocally, if \eqref{eq-rm} holds, then $\alpha$ is contained in a totally umbilical hypersurface of $\m^{n+1}(c)$.
\end{thm} 
 
 Identity \eqref{eq-rm} is valid for curves included in totally umbilical hypersurfaces $\Sigma$ of {\it any} manifold $M^{n+1}$ where, up to a constant in \eqref{eq-rm}, it holds $\sum_{i=1}^na_i\kappa_i(s)-H(\alpha(s))=0$ and $H$ is the mean curvature of $\Sigma$: see \cite[Thm. 1]{dd20}.  Thus the novelty here is that we prove the converse when $M^{n+1}$ is a space form. 
 
 Proof of these results will be done in Sect. \ref{s2}. After each one of the proofs, we will compare the techniques and the results employed with that of the existing literature to show what is actually new in this paper. In the last Section \ref{s3} we investigate those curves of $\h^3$ and $\s^3$ contained in   totally umbilical surfaces with   constant   torsion. We will derive in Thm. \ref{tx} the expression of the curvature and some numerical pictures of these curves are shown.

\section{Proof of the results}\label{s2}

 Let $\alpha\colon I\subset\r\to \m^3(c)$ be a regular (non-geodesic) curve  parametrized by arc length.  Let $\kappa>0$ and $\tau$ be the curvature and the torsion of $\alpha$, respectively. Associated to $\alpha$ there is a Frenet frame $\{T,N,B\}$  whose Frenet equations are \begin{equation*}
\begin{split}
\nabla_TT&=\kappa N,\\
\nabla_TN&=-\kappa T+\tau B,\\
\nabla_TB&=-\tau N.
\end{split}
\end{equation*}
Here $\nabla$ denotes the Levi-Civita connection of $\m^3(c)$. If $\nabla^0$ is the Levi-Civita connection of $\r^4$ with respect to the   metric $\langle,\rangle$, then  $\nabla$ is $\nabla^0$ are related by the identity 
\begin{equation}\label{nn}
\nabla^{0}_{X}Y=\nabla_X Y-c\, \langle X,Y\rangle p,
\end{equation}
where $X,Y\in\mathfrak{X}(\m^3(c))$ are tangent vectors at $p\in\m^3(c)$. 
\subsection{Proof of Theorem \ref{t1}}

{\it Necessity}. Since $\alpha$ is contained in $U_{\vec{a},\sigma}^h$, then $\langle\alpha(s),\vec{a}\rangle=\sigma$ for all $s\in I$. Differentiating with respect to $s$, we have $\langle T,\vec{a}\rangle=0$. Differentiating again and using the Frenet equations together \eqref{nn}, we have 
$$\langle \nabla_TT+\langle T,T\rangle\alpha,\vec{a}\rangle=0,$$
that is
\begin{equation}\label{kn}
\kappa\langle N,\vec{a}\rangle+\sigma=0.
\end{equation}
We distinguish two cases.
\begin{enumerate}
\item Case $\sigma=0$. This corresponds when $U_{\vec{a},0}$ is a totally geodesic plane, in particular, $\epsilon=1$.  Since $\kappa\not=0$, then \eqref{kn} yields $\langle N,\vec{a}\rangle=0$ and differentiating with respect to $s$ we have $\tau\langle B,\vec{a}\rangle=0$. This implies   $\tau=0$ identically in $I$ because if $\langle B,\vec{a}\rangle=0$, using that  $\{T,N,B,\alpha\}$ is an orthonormal basis of $(\r^4,\langle,\rangle)$, we would have $1=\langle\vec{a},\vec{a}\rangle=-\langle\alpha,\vec{a}\rangle^2=0$, a contradiction.
\item Case $\sigma\not=0$. Differentiating \eqref{kn}   together the second Frenet equation and the fact that $\langle T,\vec{a}\rangle=0$, we get 
\begin{equation}\label{bb}
\kappa'\langle N,\vec{a}\rangle+\kappa\tau\langle B,\vec{a}\rangle=0.
\end{equation}
On the other hand, because $\langle\vec{a},\vec{a}\rangle=\epsilon$ and  the fact that $\{T,N,B,\alpha\}$ is an orthonormal basis of $(\r^4,\langle,\rangle)$,  we deduce  from identities \eqref{kn} and \eqref{bb}  
\begin{equation*}
\begin{split}
\epsilon&=\langle\vec{a},\vec{a}\rangle=\langle N,\vec{a}\rangle^2+\langle B,\vec{a}\rangle^2-\langle \alpha,\vec{a}\rangle^2\\
&=\frac{\sigma^2}{\kappa^2}+\frac{\sigma^2}{\kappa ^2}\frac{\kappa'^2}{\kappa^2\tau^2}-\sigma^2.
\end{split}
\end{equation*}
This implies 
$$\frac{1}{\kappa^2}+\frac{\kappa'^2}{\kappa^4\tau^2}=1+\frac{\epsilon}{\sigma^2}.$$
This identity together \eqref{hh} is just \eqref{eq1}.
\end{enumerate}

{\it Sufficiency}. We distinguish  two cases according to $\tau$.
\begin{enumerate}
\item Case $\tau=0$. Then $\nabla_TB=0$, hence $\nabla_T^0B=0$ by \eqref{nn}. This proves that $B=B(s)$ is a constant function. Let $\vec{a}=B$, in particular, $\langle\vec{a},\vec{a}\rangle=\langle B,B\rangle=1$. Since the derivative of the function $s\mapsto\langle\alpha(s),\vec{a}\rangle$ is $\langle T,\vec{a}\rangle=0$, then the function $\langle \alpha(s),\vec{a}\rangle$ is constant. Because $\vec{a}=B$ is a tangent vector of $\h^3$, then $\langle \alpha(s),\vec{a}\rangle=0$ for all $s\in I$, proving that $\alpha$ is contained   in the totally geodesic plane $U_{\vec{a},0}$. 
\item Case $\tau\not=0$. Define the function $\beta\colon I\to\r^4$ by
\begin{equation}\label{beta}
\beta(s)=\alpha(s)+\frac{1}{\kappa(s)}N(s)-\frac{\kappa'(s)}{\kappa(s)^2\tau(s)}B(s).
\end{equation}
Using   the Frenet equations, we have  
\begin{equation}\label{k2}
\frac{d}{ds}\beta(s)=\left(\frac{\tau}{\kappa}-\left(\frac{\kappa'}{\kappa^2\tau}\right)'\right)B.
\end{equation}
The derivative of \eqref{eq13} yields
\begin{equation}\label{k3}
\frac{\kappa'}{\tau\kappa^2}\left( \left(\frac{\kappa'}{\kappa^2\tau}\right)'-\frac{\tau}{\kappa}\right)=0.
\end{equation}
Since $\kappa'=0$, the parenthesis must vanish, hence \eqref{k2} leads to $\frac{d}{ds}\beta(s)=0$ for all $s\in I$. This proves that   $\beta$ is a constant function. Let $\beta(s)=\vec{b}$ for some vector $\vec{b}\in\r^4$.  
Multiplying  \eqref{beta} by $\alpha$ and because $\alpha$ is orthogonal to $N$ and $B$, we have $\langle\alpha(s),\vec{b}\rangle=-1$ for all $s\in I$.  Moreover, using \eqref{eq13} and \eqref{beta} we get
$$\langle \vec{b},\vec{b}\rangle=-1+\frac{1}{\kappa^2}+\frac{\kappa'^2}{\kappa^4\tau^2}=-1+C.$$
If $C=1$, then   $\langle \vec{b},\vec{b}\rangle=0$, and $\alpha$ is contained in the horosphere $U^h_{\vec{b},-1}$. If $C\not=1$, let $\vec{a}=\vec{b}/\sqrt{|C-1|}$. Then $\langle\vec{a},\vec{a}\rangle=1$ if $C>1$ and $\langle\vec{a},\vec{a}\rangle=-1$ if $C<1$. In both cases, $\langle\alpha(s),\vec{a}\rangle=-1/\sqrt{|C-1|}:=\sigma$. Thus $\alpha$ is contained in the equidistant surface $U^h_{\vec{a},\sigma}$  if $C>1$ or in the   geodesic sphere $U^h_{\vec{a},\sigma}$ if $C<1$. Finally, once we know that $\alpha$ is contained in a totally umbilical surface, the value of $C$ is given by \eqref{eq1}, that is, $C=1/H^2$.  
\end{enumerate}

Equation \eqref{eq1} can be also obtained from the computations done in \cite{dd20} and that the authors did not mention. Indeed, in Remark 1 of \cite{dd20} the authors employ Frenet frames for curves in a generic Riemannian $3$-manifold $M^3$. Following the calculations and notations of \cite{dd20}, let $\alpha$ be a curve contained in a totally umbilical surface $S\subset M^3$.  Let $\xi$ and $H$ be the unit  normal and the mean curvature of $S$, respectively. Since $S$ is totally umbilical, then the tangential component of $\nabla_T\xi$ vanishes. Thus $\nabla_{T}\xi=c_1 N+c_2 B$. Then $c_1=H/\kappa$, $c_2'=-\tau c_1$ and $c_1'=\tau c_2$. The same computation of \cite[Rem. 1]{dd20} gives 
$$1=\langle \xi,\xi\rangle=c_1^2+c_2^2=\left(\frac{H}{\kappa}\right)^2+\frac{1}{\tau^2}\left(\frac{H}{\kappa}\right)'^2.$$
Since $H$ is constant because $S$ is a totally umbilical surface, and assuming $H\not=0$, we obtain 
$$\frac{1}{H^2}=\left(\frac{1}{\kappa}\right)^2+\frac{1}{\tau^2}\left(\frac{1}{\kappa}\right)'^2.$$
This identity coincides with \eqref{eq1}.

\begin{rem}\label{r2} In the case that the totally umbilical surface is a   geodesic sphere of $\h^3$ or $\s^3$, part of the proof of Thm. \ref{t1} is contained in  \cite{dd}. Indeed, the vector $\t_\beta$ in Thm. 4 of \cite{dd}  is nothing but the unit normal of the sphere, and from Eqs. (27) and (30) of \cite{dd},  we can characterize spherical curves with an equation depending on a constant. This constant can be expressed in terms of   the mean curvature of the geodesic sphere.  
\end{rem}
\subsection{Proof of Theorem \ref{t12}}

We only need to consider that case $\sigma\not=0$. Necessity is immediate because \eqref{eq12} is  the derivative of Eq. \eqref{eq1} multiplied by $\kappa'/(\tau\kappa^2)$: see \eqref{k3}. But $\kappa'\not=0$ and the result follows. For sufficiency, the function $\beta$ defined in \eqref{beta} is constant because of \eqref{eq12} and \eqref{k2}. Then  the arguments given in the ``Sufficiency'' part of the proof of Thm. \ref{t1}  imply that $\alpha$ is contained in a totally umbilical surface.    

\subsection{Proof of Theorems \ref{t2} and \ref{t22}}

The proofs are similar to that of Thm. \ref{t1}  and \ref{t12}. We only show the differences in the proof of Thm. \ref{t2}.

{\it Necessity}. We have  $\langle\alpha(s),\vec{a}\rangle=\sigma$ and $\langle T,\vec{a}\rangle=0$. Differentiating twice again, we have 
$$\kappa\langle N,\vec{a}\rangle+\sigma=0$$
and
$$\kappa'\langle N,\vec{a}\rangle+\kappa\tau\langle B,\vec{a}\rangle=0.$$
The case $\sigma=0$ follows as in Thm. \ref{t1}. If $\sigma\not=0$ and because $\tau\not=0$, we have
Thus
\begin{equation*}
\begin{split}
1&=\langle\vec{a},\vec{a}\rangle=\langle N,\vec{a}\rangle^2+\langle B,\vec{a}\rangle^2+\langle \alpha,\vec{a}\rangle^2 =\frac{\sigma^2}{\kappa^2}+ \sigma^2 \frac{\kappa'^2}{\kappa^4\tau^2}+\sigma^2\\
&=\sigma^2\left(\frac{1}{\kappa^2}+\frac{\kappa'^2}{\kappa^4\tau^2}+1\right).
\end{split}
\end{equation*}
This gives \eqref{eq2}.

{\it Sufficiency}. The case $\tau=0$ follows with the same steps. If $\tau\not=0$, differentiating \eqref{eq2} we obtain \eqref{eq3} which allows to see that the function $\beta$ defined in \eqref{beta} is constant, $\beta(s)=\vec{b}$. Now it is immediate from \eqref{eq14} that $|\vec{b}|^2=1+C$ and $\langle \alpha(s),\vec{b}\rangle=1$ for all $s\in I$. Letting $\vec{a}=\frac{\vec{b}}{\sqrt{1+C}}$, we obtain $\langle\vec{a},\vec{a}\rangle=1$ and $\langle\alpha(s),\vec{a}\rangle=1/\sqrt{1+C}$. This proves that  $\alpha$ is contained in the geodesic sphere $U^s_{\vec{b},\sigma}$ where $\sigma=1/\sqrt{1+C}$. In particular, $C=\frac{1-\sigma^2}{\sigma^2}$ which coincides with $1/H^2$, where $H$ is the mean curvature of $U^s_{\vec{b},\sigma}$.

\subsection{Proof of Theorems \ref{t5}}

First we recall that RM frames along regular curves are frames of type $\{T,N_1,\ldots,N_n\}$, where $T$ and $\nabla_TN_i$ are parallel. Therefore, the RM equations are 
\begin{equation}\label{eqrm}
\begin{split}
\nabla_TT&=\sum_{i=1}^n\kappa_i N_i,\\
\nabla_TN_i&=-\kappa_i N_i,\quad 1\leq i\leq n.
\end{split}
\end{equation}
The functions $\kappa_i$ are called the curvatures associated to the RM frame. 

As we pointed out in the Introduction, the first part of the statement    was done in \cite{dd20}. For the sake of completeness, we recall it. Since $\langle\alpha(s),\vec{a}\rangle=\sigma$, the derivative with respect to $s$ gives $\langle T,\vec{a}\rangle=0$. Differentiating with respect to $s$ again and using \eqref{nn} and \eqref{eqrm}, we have 
\begin{equation}\label{eq9}
\sum_{i=1}^n\kappa_i\langle N_i,\vec{a}\rangle+\sigma=0.
\end{equation}
On the other hand, from \eqref{eqrm} one gets
$$\frac{d}{ds}\langle N_i,\vec{a}\rangle=\langle\nabla_TN_i,\vec{a}\rangle=-\kappa_i\langle T,\vec{a}\rangle=0.$$
Thus $\langle N_i,\vec{a}\rangle$ is a constant function $a_i$ for all $1\leq i\leq n$. This together \eqref{eq9} proves \eqref{eq-rm}.

We now prove the converse. Define $\beta(s)=\sum_{i=1}^na_iN_i(s)-\sigma\alpha(s)$. Then \eqref{eqrm} implies
\begin{equation*}
\begin{split}
\beta'(s)&=\sum_{i=1}^na_i\nabla_T N_i(s)-\sigma T(s)=-\sum_{i=1}^na_i\kappa_i(s)T(s)-\sigma T(s)\\
&=0.
\end{split}
\end{equation*}
Thus there is $\vec{a}\in \r^{n+2}=\m^{n+1}(c)$ such that $\beta(s)=\vec{a}$ for all $s\in I$, that is, 
$$\sum_{i=1}^na_iN_i(s)-\sigma\alpha(s)=\vec{a}.$$
In particular, if we multiply by $\alpha(s)$ we obtain $ \langle\alpha(s),\vec{a}\rangle=-c\sigma$. This proves that $\alpha$ is contained in the totally umbilical hypersurface $U_{\vec{a},-c\sigma}$.  

\begin{rem}Theorem \ref{t5} shows how powerful RM  frames are in comparison with the Frenet frame. Indeed, from Remark 2 of \cite{dd} we have a taste of how complicated characterizing spherical curves in higher dimensions can get. For curves in $\r^{n+1}$, it is possible to find a generic equation characterizing spherical curves using in Frenet frames following Thm. 5.1 of \cite{dsf21} which generalizes the given one in Rem. 2 of \cite{dd20}.
\end{rem}

\begin{rem} Equation \eqref{eq-rm} can be interpreted as model-dependent because of $\sigma$ in the right hand-side. However, we can renormalize the linear relation \eqref{eq-rm} such that $a_1^2+\ldots+a_n^2 = 1$. In such a case,  the constants $a_i$  can be interpreted as the coordinates of the unit normal of the hypersurface along the curve whereas $\sigma$ changes to be the mean curvature of the totally umbilical hypersurface. Notice that in \cite[Thm. 3]{dd}, the value of $\coth(R)$ (resp. $\cot(R)$) is the mean curvature of the geodesic hypersphere of radius $R$ in $\h^{n+1}$ (resp. $\s^{n+1}$): see also  \cite[Thm. 1]{dd20}. 
\end{rem}

\section{Curves  with constant torsion contained in totally umbilical surfaces}\label{s3}
In this section we study curves in  a totally umbilical surface of $\h^3$ and $\s^3$ with   constant   torsion. If instead of curves with constant torsion   we ask for those curves with constant curvature $\kappa$, then the curves are known. If the curve is included in a totally umbilical surface $U_{\vec{a},\sigma}$ with constant mean curvature $H$ then $\kappa=H+\kappa_g$, where $\kappa_g$ is the geodesic curvature of the curve. Thus we deduce that $\kappa$ is constant if and only if  $\kappa_g$ is constant. The curves of constant geodesic curvature in $U_{\vec{a},\sigma}$ are known because $U_{\vec{a},\sigma}$ is isometric to a space form. If  the ambient space is $\h^3$, then it is isometric to $\r^2$,  a hyperbolic plane $\h^2(c)$ and a sphere $\s^2(c)$  in case of a  horosphere, geodesic plane or equidistant surface,  and geodesic sphere, respectively. If $U_{\vec{a},\sigma}$ is included in $\s^3$, then it is isometric to a sphere $\s^2(c)$.

Suppose now that $\alpha$ is contained in a totally umbilical surface $U_{\vec{a},\sigma} $ of $\m^3(c)$ and that the  torsion $\tau$ of $\alpha$  is constant. When the ambient space is the Euclidean space $\r^3$,   spherical curves with   constant torsion are studied in \cite{ks} and  parametrizations of such curves are obtained in terms of hypergeometric functions. We  extend these results in the context of totally umbilical surfaces of space forms. We will discard the case $\tau=0$ identically because then $\alpha$ is included in   a totally geodesic surface and in such a case, the curvature is an arbitrary function.   Assume the case that   the torsion $\tau$  is a non-zero constant and  we will derive explicitly the curvature $\kappa$ of the curve. We will exclude the particular case that $\kappa$ is constant which it has been previously considered.

 \begin{thm}\label{tx} Let $U_{\vec{a},\sigma} $ be a totally umbilical surface of $\h^3$ or $\s^3$. If $\alpha$   is a regular curve  contained in $U_{\vec{a},\sigma}$ with constant non-zero torsion $\tau$ and $\kappa'\not=0$, then the curvature of $\alpha$ is 
 \begin{equation}\label{sol1}
 \kappa(s)=\frac{|H|}{ \sin(\tau s+a)},\quad a\in\r,
 \end{equation} 
 where $H$ is the mean curvature of $U_{\vec{a},\sigma} $.
 \end{thm}

It will always be understood that the domain of $\alpha$ is maximal in order to have well defined functions  and that the curvature  $\kappa$ is positive. 

\begin{proof}
Let $\tau\not=0$ be a real number. Since the torsion is non-zero, then $U_{\vec{a},\sigma} $ is not a totally geodesic surface, in particular, $H\not=0$. Then Eqs. \eqref{eq1} and \eqref{eq2} can be written as  
$$\tau^2\left(\frac{1}{\kappa}\right)^2+{\left(\frac{1}{\kappa}\right)'}^2=\frac{\tau^2}{H^2}.$$
Let $w(s)=1/\kappa(s)$. Notice that $w$ is not a constant function because $\kappa'\not=0$. Then the above  equation reads as $\tau^2 w^2+w'^2=\tau^2/H^2$ whose solution is 
$w(s)=\frac{1}{|H|}\sin(\tau s+a)$, $a\in\r$. This gives  \eqref{sol1}.
\end{proof}

  Once we know $\kappa$ and $\tau$, the curve $\alpha$ is determined up to rigid motions of $\m^3(c)$. A further step would be to get an explicit parametrization of the curve $\alpha$. This will involve the solution of \eqref{sol1} together the fact that $\alpha$ is parametrized by arc-length. Such as it occurs   in   \cite{ks}, hypergeometric functions are used. However, our interest in this section is  to give numerical pictures of curves with   constant torsion contained in totally umbilical surfaces  in $\h^3$ and $\s^3$. We will employ the software Mathematica to make these figures \cite{mat}.

We begin with curves in hyperbolic space $\h^3$. In order to visualize the curve, we will show the picture in the upper-halfspace model $\r_+^3$ of $\h^3$ where  $\alpha$ will be a curve in the $3$-dimensional vector space $\r^3_+$ instead of $\mathcal{H}$.  The numerical  computations will be done first in the hyperquadric model $\mathcal{H}$ and finally, the isometry
$$\Phi\colon \mathcal{H} \to\r^3_+,\quad \Phi(x_1,x_2,x_3,x_4)=\left(\frac{x_1}{x_3+x_4},\frac{x_2}{x_3+x_4},\frac{1}{x_3+x_4}\right)$$
between both models of $\h^3$ will give us the desired picture. 

\begin{figure}[hbtp]
\centering
\includegraphics[width=.75\textwidth]{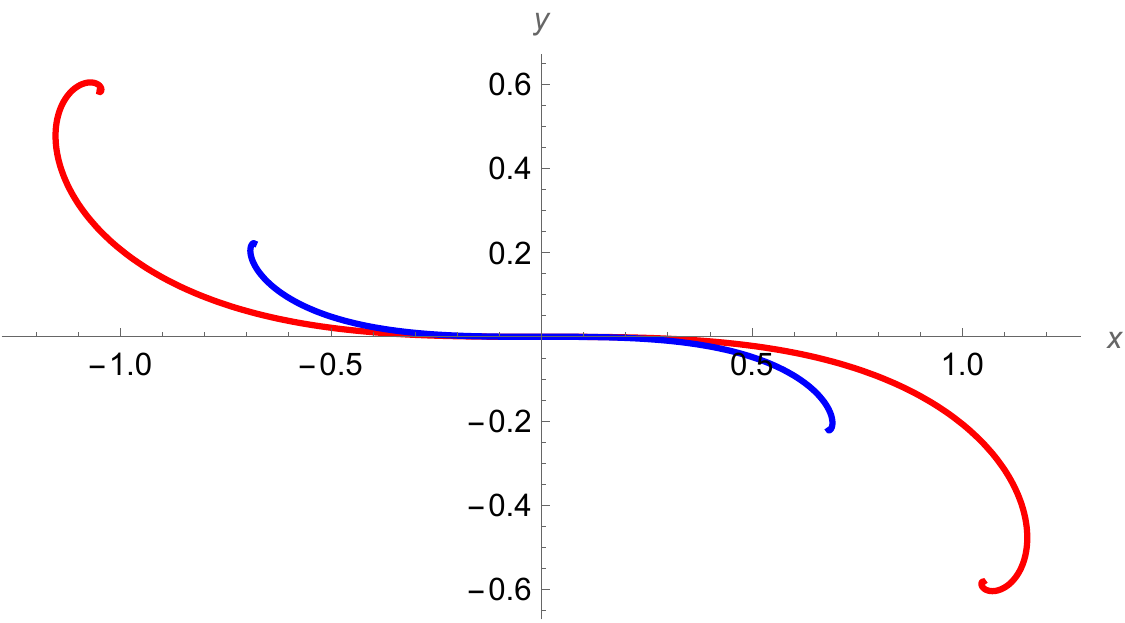}
\caption{Two curves contained in the horosphere $z=1$ in the upper-halfspace model $\r_+^3$ of $\h^3$ with constant torsion $\tau=1$ (red) and $\tau=2$ (blue). The initial conditions are $x_1(s_0)=x_2(s_0)=0$ where $s_0=\pi/2$ ($\tau=1)$ and $s_0=\pi/4$ ($\tau=2$).}\label{fig1}
\end{figure}

In the example that we show in this section,   we take a horosphere of $\h^3$ as the totally umbilical surface. After an isometry of $\h^3$, a horosphere in $\r_+^3$ is a horizontal plane. Thus  $\alpha$ will be a curve contained in a (horizontal) plane of $\r^3$ which we identify with $\r^2$. Consider in $\mathcal{H}$ the horosphere $U^h_{\vec{a},1}$, where $\vec{a}=(0,0,1,-1)$. The equations of $U^h_{\vec{a},1}$ are $x_1^2+x_2^2+x_3^2-x_4^2=1$ and $x_3+x_4=1$. This horosphere, via $\Phi$, is the horizontal plane  of equation $z=1$ in $\r^3_+$, where $(x,y,z)$ are the standard coordinates of $\r^3_+$. 
Let $\tau$ be the a non-zero real number which will be the torsion of the curve that we are looking for. Since $|H|=1$, then \eqref{sol1} gives $\kappa(s)=\frac{1}{\sin (\tau s+a)}$. Take $a=0$. Let 
\begin{equation}\label{c1}
\alpha(s)=(x_1(s),x_2(s),x_3(s),1-x_3(s))
\end{equation}
 be a curve in $\h^3= \mathcal{H}$ parametrized by arc-length with torsion $\tau$. This gives two equations, namely, 
\begin{equation*}
\begin{split}
x_1(s)^2+x_2(s)^2+2x_3(s)&=0,\quad  \mbox{because }\alpha(s)\in\h^3 ,\\
x_1'(s)^2+x_2'(s)^2&=1,\quad  \mbox{because } |\alpha'(s)|=1 .
\end{split}
\end{equation*}
 From the second equation, there is a smooth function $\theta=\theta(s)$ such that 
 \begin{equation}\label{bb2}
 \left\{
\begin{split}
x_1'(s) &=\cos\theta(s),\\
x_2'(s) &=\sin\theta(s).
\end{split}\right.
\end{equation}
We now calculate the curvature of $\alpha$. We know that $\kappa(s)=|\nabla_{\alpha'(s)}\alpha'(s)|$. From \eqref{nn}, we have 
\begin{equation*}
\begin{split}
\nabla_{\alpha'(s)}\alpha'(s)&=\nabla^0_{\alpha'(s)}\alpha'(s)-\langle\alpha'(s),\alpha'(s)\rangle\alpha(s)=\alpha''(s)-\alpha(s)\\
&=(x_1''(s)-x_1(s),x_2''(s)-x_2(s),x_3''(s)-x_3(s),-x_3''(s)+x_3(s)-1).
\end{split}
\end{equation*}
Using both equations \eqref{bb2}, we obtain 
\begin{equation}\label{kh}
\kappa(s)=|\nabla_{\alpha'(s)}\alpha'(s)|=\sqrt{1+\theta'(s)^2}.
\end{equation}
Knowing that $\kappa(s)=1/\sin (\tau s)$, the above equation reduces into $\theta'(s)=\pm \cot (\tau s)$. The solution of this equation is  
$\theta(s)=\frac{1}{\tau}\log\sin (\tau s)$ after choices of integration constants.  Finally the curve $\alpha$ given in \eqref{c1}  is determined by the ODE system
 \begin{equation*}
\begin{split}
x_1'(s) &=\cos(\frac{1}{\tau}\log \sin (\tau s)),\\
x_2'(s) &=\sin(\frac{1}{\tau}\log\sin (\tau s)).
\end{split}
\end{equation*}
Once we know $x_1(s)$ and $x_2(s)$, the curve $\alpha$ given in \eqref{c1} is completely determined by the identities 
$$x_3(s)=-\frac{x_1(s)^2+x_2(s)^2}{2},\quad x_4(s)=1+\frac{x_1(s)^2+x_2(s)^2}{2}.$$
 We show  in Fig. \ref{fig1} two curves of constant torsions $\tau=1$ and $\tau=2$ for given initial conditions. Here we are identifying the horosphere $z=1$ with $\r^2$. Since $\kappa(s)=1/\sin(\tau s)$, the domain of the curve is $(0,\frac{\pi}{\tau})$. In particular, as $s\to 0$ and $s\to\frac{\pi}{\tau}$, we have $\kappa\to\infty$. This explains that the curve in Fig. \ref{fig1} spirals  very fast close to the end points of the curve.

We now present an example in $\s^3$. For the geodesic sphere $U^s_{\vec{a},\sigma}$, let $\vec{a}=(0,0,0,1)$ and $\sigma=1/2$, in particular, $|H|=1/\sqrt{3}$. Then 
$$U^s_{\vec{a},\sigma}=\{(x_1,x_2,x_3,\frac12)\in\s^3:x_1^2+x_2^2+x_3^2=1-\sigma^2=\frac{3}{4}\}.$$
 In order to visualize this surface, we identify $U^s_{\vec{a},\sigma}$ with the sphere 
$\s^2(R)=\{(x,y,z)\in\r^3\colon x^2+y^2+z^2=R^2\}$ of radius $R=\sqrt{3}/2$ in Euclidean space $\r^3$ thanks to the map
$$(x_1,x_2,x_3,\frac12)\in U^s_{\vec{a},\sigma}\longmapsto (x_1,x_2,x_3)\in \s^2(R).$$
 Again, let $\tau$ be a non-zero real number. By choosing $a=0$, then the expression  \eqref{sol1} gives $\kappa(s)=\frac{1}{\sqrt{3}\sin (\tau s)}$, where we are assuming $\sin(\tau s)>0$ in order to ensure that $\kappa$ is positive.   

\begin{figure}[hbtp]
\centering
\includegraphics[width=.5\textwidth]{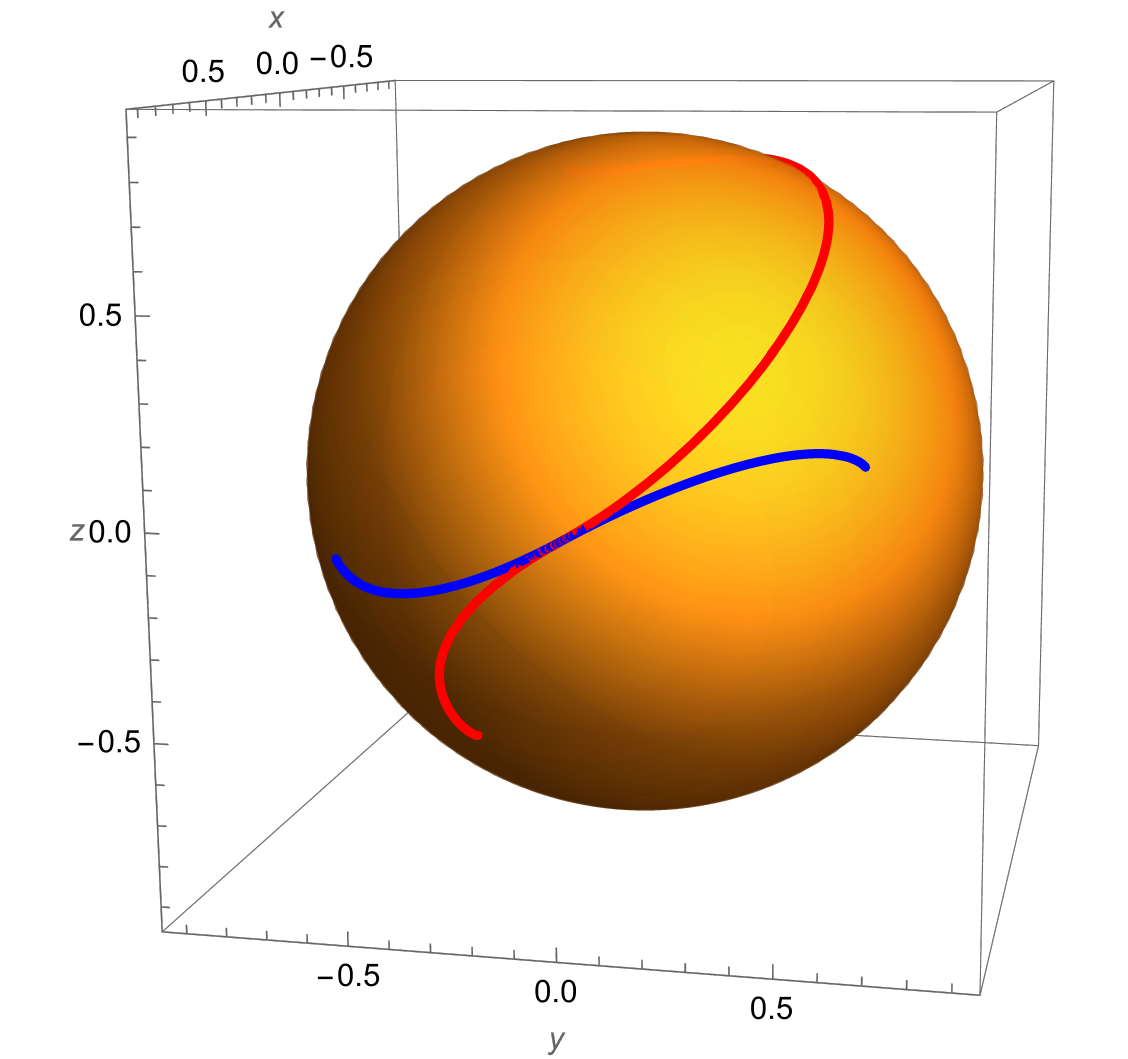}
\caption{Two curves with constant torsion contained in the sphere $\s^2(R)$, $R=\sqrt{3}/{2}$, identifying $\s^2(R)$ with $U^s_{\vec{a},1/2}$.  The  torsion is $\tau=1$ (red) and $\tau=2$ (blue). The initial conditions are given at $s_0=\pi/4$:   $\varphi(s_0)=\theta(s_0)=0$, $\theta'(s_0)=1$ and $\varphi'(s_0)=\frac{1}{\sqrt{3}}$. }\label{fig2}
\end{figure}

 Let $\alpha(s)=(x(s),y(s),z(s))$ be a curve in $\s^2(R)$ parametrized by arc-length. Since $\alpha(s)\in\s^2(R)$ for all $s\in I$, we parametrize $\alpha$ by  
\begin{equation}\label{as}
\alpha(s)=R(\cos\varphi(s)\cos\theta(s), \cos\varphi(s)\sin\theta(s),\sin\varphi(s)).
\end{equation}
The fact that $\alpha$ is parametrized by arc-length writes as 
\begin{equation}\label{ss2}
\varphi'^2+\theta'^2(\cos\varphi)^2=\frac{4}{3}.
\end{equation}
For the computation of $\kappa$ we know $\kappa=|\nabla_TT|$.   From \eqref{nn}, we have 
\begin{equation*}
\nabla_{\alpha'(s)}\alpha'(s)=\nabla^0_{\alpha'(s)}\alpha'(s)+\langle\alpha'(s),\alpha'(s)\rangle\alpha(s)=\alpha''(s)+\alpha(s).
\end{equation*}
Then equation $\kappa=|\nabla_TT|$ becomes
\begin{equation}\label{ks2}
\begin{split}
 &(\cos \varphi )^2\theta''^2+\varphi ''^2-2 \theta ' \theta '' \varphi ' \sin (2 \varphi )+\theta '^2 \sin(2\varphi)\varphi ''\\
 &-\theta '^2 \left(\varphi '^2 (\cos (2 \varphi )-3)+\cos (2 \varphi )+1\right) +\theta '^4 (\cos \varphi )^2+\varphi '^4-2 \varphi '^2+1\\
 &=\frac{4}{9(\sin (\tau s))^2}.
 \end{split}
  \end{equation}
This differential equation of second order, together \eqref{ss2}, are numerically solved after initial conditions are given for the functions $\varphi$ and $\theta$ and their first derivatives at $s=s_0$.  In  Fig. \ref{fig2} we show two  curves with torsion $\tau=1$ and $\tau=2$. In both cases, we take $s_0=\pi/4$ and the domain for the variable $s$ is $(0,\pi)$ and $(0,\frac{\pi}{2})$ respectively. As in the case of Fig. \ref{fig1}, at the end points of the domain the curvature of the curve goes to $\infty$ and thus, the curve spirals very fast as well as the parameter to the boundary points of the intervals $(0,\pi)$ and $(0,\frac{\pi}{2})$, respectively.

\begin{rem} 
\begin{enumerate}
\item The curve $\alpha$  given in \eqref{as} and  satisfying \eqref{ks2}   has constant torsion when $\alpha$  is   considered as a curve in the ambient space $\s^3$ after the identification with the sphere $\s^3\cap\{x_4=\frac12\}$. In particular, $\alpha$ is not  one of the curves   given in \cite{ks}.
\item Comparing   the computations of $\s^3$ with  the example of the curve included in a horosphere of $\h^3$ (Fig. \ref{fig1}),     the calculation of Eq. \eqref{ks2} is  different. Notice  that simplicity of \eqref{kh} is in part due to the horosphere $z=1$ is isometric to the Euclidean plane $\r^2$.  
\end{enumerate}

\end{rem}

\section*{Acknowledgements}
 The author would like to thank the anonymous referee for insightful suggestions, which helped to improve the exposition of the paper. In particular, for pointing out the results of papers \cite{dd,dd20}. 
 
The author  is a member of the IMAG and of the Research Group ``Problemas variacionales en geometr\'{\i}a'',  Junta de Andaluc\'{\i}a (FQM 325). This research has been partially supported by MINECO/MICINN/FEDER grant no. PID2023-150727NB-I00,  and by the ``Mar\'{\i}a de Maeztu'' Excellence Unit IMAG, reference CEX2020-001105- M, funded by MCINN/AEI/10.13039/ 501100011033/ CEX2020-001105-M.

\section*{Declarations statements}

{\bf Competing interests.} The author declares no competing interests.

\noindent {\bf Data availability.} No datasets were generated or analyzed during the current study.

\end{document}